\documentclass[a4paper,12pt,oneside]{article}
\usepackage{amsmath,amsfonts,amssymb,amsmath,latexsym,amsthm,enumerate}

\newtheorem{thm}{Theorem}

\theoremstyle{definition}

\theoremstyle{plain}
\newtheorem{lem}[thm]{Lemma}

\newtheorem{prop}[thm]{Proposition}

\begin{document}
\title {$q$-Bernoulli polynomials and $q$-umbral calculus}
\author{by \\Dae San Kim and Taekyun Kim}\date{}\maketitle

\begin{abstract}
\noindent In this paper, we investigate some properties of $q$-Bernoulli polynomials arising from $q$-umbral calculus. Finally, we derive some interesting identities of $q$-Bernoulli polynomials from our investigation.
\end{abstract}

\section{Introduction and preliminaries}

Throughout this paper we will assume $q$ to be a fixed number between $0$ and $1$.  We denote by $D_{q}$ the $q$-derivative of a function
\begin{equation}\label{eq:1}
\left(D_{q}f\right)(x)=\frac{f\left(qx\right)-f\left(x\right)}{\left(q-1\right)x},\,\,\,(\text{see}\,\, \lbrack 8,10\rbrack).
\end{equation}
The Jackson definite $q$-integral of the function $f$ is defined by
\begin{equation}\label{eq:2}
\int_{0}^{x}f(t)d_{q}t=\left(1-q\right)\sum_{a=0}^{\infty}f\left(q^{a}x\right)xq^{a},\,\,\,(\text{see}\,\, \lbrack 8,12,13\rbrack).
\end{equation}
From (\ref{eq:1}) and (\ref{eq:2}), we note that
\begin{equation*}
D_{q}\int_{0}^{x}f(t)d_{q}t=f(x),\quad \int_{a}^{b}f(x)d_{q}x=\int_{0}^{b}f(x)d_{q}x-\int_{0}^{a}f(x)d_{q}x.
\end{equation*}
In this paper, we use the following notations:
\begin{equation}\label{eq:3}
\lbrack x\rbrack_{q}=\frac{1-q^{x}}{1-q},\quad \left(a+b\right)_{q}^{n}=\prod_{i=0}^{n-1}\left(a+q^{i}b\right),\,\,\,\,\left(n\in\mathbf{Z}_{+}\right)
\end{equation}
and
\begin{equation}\label{eq:4}
\left(1+a\right)_{q}^{\infty}=\prod_{j=0}^{\infty}\left(1+q^{j}a\right),\quad \lbrack n \rbrack_{q}!=\lbrack n\rbrack_{q}\lbrack n-1\rbrack_{q}\cdots \lbrack 2\rbrack_{q}\lbrack 1\rbrack_{q}.
\end{equation}
The $q$-analogue of exponential function is defined by
\begin{equation}\label{eq:5}
e_{q}(t)=\frac{1}{\left(1-\left(1-q\right)t\right)_{q}^{\infty}}=\sum_{n=0}^{\infty}\frac{t^{n}}{\lbrack n\rbrack_{q}!},\,\,\,(\text{see}\,\, \lbrack 5,6,8,10\rbrack).
\end{equation}
In $\lbrack 10\rbrack$, the $q$-analogues of Bernoulli polynomials are defined by the generating function to be
\begin{equation}\label{eq:6}
\frac{t}{e_{q}(t)-1}e_{q}(xt)=\sum_{n=0}^{\infty}B_{n,q}(x)\frac{t^{n}}{\lbrack n\rbrack_{q}!},\,\,\,(\text{see}\,\, \lbrack 8-14\rbrack).
\end{equation}
In the special case, $x=0$, $B_{n,q}(0)=B_{n,q}$ is called the $n$-th $q$-Bernoulli number.\\
From (\ref{eq:6}), we can derive the following equation:
\begin{equation}\label{eq:7}
B_{n,q}(x)=\sum_{l=0}^{n}\binom{n}{l}_{q}x^{n-l}B_{l,q}=\sum_{l=0}^{n}B_{n-l,q}x^{l}\binom{n}{l}_{q},
\end{equation}
where $\binom{n}{l}_{q}=\frac{\lbrack n\rbrack_{q}!}{\lbrack l\rbrack_{q}!\lbrack n-l\rbrack_{q}!}=\frac{\lbrack n\rbrack_{q}\lbrack n-1\rbrack_{q}\cdots \lbrack n-l+1\rbrack_{q}}{\lbrack l\rbrack_{q}!}.$\\
Let $\mathbf{C}$ be the complex number field and let $\mathcal{F}$ be the set of all formal power series in variable $t$ over $\mathbf{C}$ with
\begin{equation}\label{eq:8}
\mathcal{F}=\left\{f(t)=\sum_{k=0}^{\infty}\frac{a_{k}}{\lbrack k\rbrack_{q}!}t^{k}\Bigg\vert a_{k}\in\mathbf{C}\right\}.
\end{equation}
Let $\mathbb{P}=\mathbf{C}\lbrack t\rbrack$ and let $\mathbb{P}^{*}$ be the vector space of all linear functionals on $\mathbb{P}$. Now we denote by $\left\langle L\vert p(x)\right\rangle$ the action of the linear functional $L$ on the polynomial $p(x)$. We remind that the vector space operations on $\mathbb{P}^{*}$ are defined by
\begin{equation*}
\left\langle L+M\vert p(x)\right\rangle=\left\langle L\vert p(x)\right\rangle+\left\langle M\vert p(x)\right\rangle,\quad \left\langle cL\vert p(x)\right\rangle=c\left\langle L\vert p(x)\right\rangle,
\end{equation*}
where $c$ is any constant in $\mathbf{C}$ (see $\lbrack 15,16\rbrack$).\\
For $f(t)=\sum_{k=0}^{\infty}\frac{a_{k}}{\lbrack k\rbrack_{q}!}t^{k}\in\mathcal{F}$, we define the linear functional on $\mathbb{P}$ by setting
\begin{equation}\label{eq:9}
\left\langle f(t)\vert x^{n}\right\rangle=a_{n}\,\,\,\,\text{for all}\,\, n\geq 0.
\end{equation}
Thus, by (\ref{eq:8}) and (\ref{eq:9}), we note that
\begin{equation}\label{eq:10}
\left\langle t^{k}\vert x^{n}\right\rangle=\lbrack n\rbrack_{q}!\delta_{n,k},\,\,\,\,(n,k\geq 0),
\end{equation}
where $\delta_{n,k}$ is the Kronecker's symbol.\\
Let $f_{L}(t)=\sum_{k=0}^{\infty}\frac{\left\langle L\vert x^{k}\right\rangle}{\lbrack k\rbrack_{q}!}t^{k}$. Then, by (\ref{eq:8}) and (\ref{eq:9}), we see that $\left\langle f_{L}(t)\vert x^{n}\right\rangle=\left\langle L\vert x^{n}\right\rangle$ and so as linear functionals $L=f_{L}(t)$. It is easy to show that the map $L\longmapsto f_{L}(t)$ is a vector space isomorphism from $\mathbb{P}^{*}$ onto $\mathcal{F}$. Henceforth, $\mathcal{F}$ denotes both the algebra of formal power series in $t$ and the vector space of all linear functionals on $\mathbb{P}$, and so an element $f(t)$ of $\mathcal{F}$ is thought of as both a formal power series and a linear functional. We call $\mathcal{F}$ the $q$-umbral algebra. The $q$-umbral calculus is the study of $q$-umbral algebra. By (\ref{eq:5}) and (\ref{eq:10}), we easily see that $\left\langle e_{q}(yt)\vert x^{n}\right\rangle=y^{n}$ and so $\left\langle e_{q}(yt)\vert p(x)\right\rangle=p(y)$.\\
Notice that for all $f(t)$ in $\mathcal{F}$
\begin{equation}\label{eq:11}
f(t)=\sum_{k=0}^{\infty}\frac{\left\langle f(t)\big\vert x^{k}\right\rangle}{\lbrack k\rbrack_{q}!}t^{k},
\end{equation}
and for all polynomials $p(x)$
\begin{equation}\label{eq:12}
p(x)=\sum_{k=0}^{\infty}\frac{\left\langle t^{k}\big\vert p(x)\right\rangle}{\lbrack k\rbrack_{q}!}x^{k},\,\,\,(\text{see}\,\, \lbrack 15,16\rbrack).
\end{equation}
For $f_{1}(t), f_{2}(t),\cdots, f_{n}(t)\in\mathcal{F}$, we have
\begin{equation}\label{eq:13}
\left\langle f_{1}(t)\cdots f_{m}(t)\big\vert x^{n}\right\rangle=\sum_{i_{1}+\cdots +i_{m}=n}\binom{n}{i_{1},\cdots,i_{m}}_{q}\left\langle f_{1}(t)\big\vert x^{i_{1}}\right\rangle\cdots\left\langle f_{m}(t) \big\vert x^{i_{m}}\right\rangle,
\end{equation}
where $\binom{n}{i_{1},\cdots,i_{m}}_{q}=\frac{\lbrack n \rbrack_{q}!}{\lbrack i \rbrack_{q}!\cdots \lbrack i_{m} \rbrack_{q}!}$.\\
The order $O\left(f(t) \right)$ of the power series $f(t)(\neq0)$ is the smallest integer $k$ for which $a_{k}$ does not vanish. If $O(f(t))=0$, then $f(t)$ is called an invertible series. If $O(f(t))=1$, then $f(t)$ is called a delta series.\\
Let $p^{(k)}(x)=D_{q}^{k}p(x)$. Then, by (\ref{eq:12}), we get
\begin{equation}\label{eq:14}
p^{(k)}(x)=\sum_{l=k}^{\infty}\frac{\left\langle t^{l}\big\vert p(x)\right\rangle}{\lbrack l\rbrack_{q}!}\lbrack l \rbrack_{q}\lbrack l-1 \rbrack_{q}\cdots \lbrack l-k+1 \rbrack_{q}x^{l-k}.
\end{equation}
From (\ref{eq:14}), we have
\begin{equation}\label{eq:15}
p^{(k)}(0)=\left\langle t^{k}\big\vert p(x)\right\rangle\,\, \text{and}\,\, \left\langle 1\big\vert p^{(k)}(x)\right\rangle=p^{(k)}(0).
\end{equation}
By (\ref{eq:15}), we get
\begin{equation}\label{eq:16}
t^{k}p(x)=p^{(k)}(x)=D_{q}^{k}p(x).
\end{equation}
Let $f(t), g(t) \in\mathcal{F}$ with $O\left(f(t)\right)=1$ and $O\left(g(t)\right)=0$. Then there exists a unique sequence $s_{n}(x)$ $\left(\deg{s_{n}(x)}=n\right)$ of polynomials such that $\left\langle g(t)f(t)^{k}\vert s_{n}(x)\right\rangle=\lbrack n\rbrack_{q}!\delta_{n,k}$, ($n,k\geq 0$), which is denoted by $s_{n}(x)\sim\left(g(t),f(t)\right)$. The sequence $s_{n}(x)$ is called the $q$-Sheffer sequence for $\left(g(t),f(t)\right)$. For $h(t), f(t), g(t)\in\mathcal{F}$ and $p(x)\in\mathbb{P}$, we have
\begin{equation}\label{eq:17}
h(t)=\sum_{k=0}^{\infty}\frac{\left\langle h(t)\vert s_{k}(x)\right\rangle}{\lbrack k\rbrack_{q}!}g(t)f(t)^{k},\quad p(x)=\sum_{k=0}^{\infty}\frac{\left\langle g(t)f(t)^{k}\vert p(x)\right\rangle}{\lbrack k\rbrack_{q}!}s_{k}(x),
\end{equation}
and
\begin{equation}\label{eq:18}
\frac{1}{g\left(\bar{f}(t)\right)}e_{q}\left(y\bar{f}(t)\right)=\sum_{k=0}^{\infty}\frac{s_{k}(y)}{\lbrack k\rbrack_{q}!}t^{k},\,\, \text{for all}\,\, y\in\mathbf{C},
\end{equation}
where $\bar{f}(t)$ is the compositional inverse of $f(t)$, (see $\lbrack 15,16\rbrack$).\\
Recently, several authors have studied $q$-Bernoulli and Euler polynomials (see $\lbrack 1-17\rbrack$). In this paper, we investigate some properties of $q$-Bernoulli polynomials arising from $q$-umbral calculus. Finally, we derive some interesting identities of $q$-Bernoulli polynomials from our results.


\section{$q$-Bernoulli polynomials and $q$-umbral calculus}

From (\ref{eq:6}), we note that
\begin{equation}\label{eq:19}
B_{n,q}(x)\sim\left(\frac{e_{q}(t)-1}{t}, t\right).
\end{equation}
By (\ref{eq:19}), we get
\begin{equation}\label{eq:20}
B_{n,q}(x)=\left(\frac{t}{e_{q}(t)-1}\right)x^{n},\,\,(n\geq 0).
\end{equation}
From (\ref{eq:7}) and (\ref{eq:16}), we note that
\begin{equation}\label{eq:21}
tB_{n,q}(x)=D_{q}B_{n,q}(x)=\lbrack n\rbrack_{q}B_{n-1,q}(x).
\end{equation}
By (\ref{eq:1}) and (\ref{eq:10}), we easily see that
\begin{align}\label{eq:22}
\left\langle\frac{e_{q}(t)-1}{t}\Bigg\vert x^{n}\right\rangle&=\frac{1}{\lbrack n+1\rbrack_{q}}\left\langle\frac{e_{q}(t)-1}{t}\Bigg\vert tx^{n+1}\right\rangle\\
&=\frac{1}{\lbrack n+1\rbrack_{q}}\left\langle e_{q}(t)-1\big\vert x^{n+1}\right\rangle=\frac{1}{\lbrack n+1\rbrack_{q}}\nonumber\\
&=\int_{0}^{1}x^{n}d_{q}x.\nonumber
\end{align}
Thus, from (\ref{eq:22}), we have
\begin{equation}\label{eq:23}
\left\langle\frac{e_{q}(t)-1}{t}\Bigg\vert p(x)\right\rangle=\int_{0}^{1}p(x)d_{q}x,\,\,\,\text{for}\,\,p(x)\in\mathbb{P}.
\end{equation}
In particular, if we take $p(x)=B_{n,q}(x)$, then
\begin{align}\label{eq:24}
\int_{0}^{1}B_{n,q}(x)d_{q}x&=\left\langle\frac{e_{q}(t)-1}{t}\Bigg\vert B_{n,q}(x)\right\rangle=\left\langle 1\Bigg\vert\frac{e_{q}(t)-1}{t}B_{n,q}(x)\right\rangle\\
&=\left\langle t^{0}\big\vert x^{n}\right\rangle=\lbrack n\rbrack_{q}!\delta_{n,0}.\nonumber
\end{align}
From (\ref{eq:7}), we can derive
\begin{align}\label{eq:25}
\int_{0}^{1}B_{n,q}(x)d_{q}x&=\sum_{k=0}^{n}B_{n-k,q}\binom{n}{k}_{q}\int_{0}^{1}x^{k}d_{q}x\\
&=\sum_{k=0}^{n}\frac{B_{n-k,q}}{\lbrack k+1\rbrack_{q}}\binom{n}{k}_{q}\nonumber
\end{align}
Therefore, by (\ref{eq:24}) and (\ref{eq:25}), we obtain the following proposition.
\begin{prop}\label{eq:prop1}
For $n\in\mathbf{Z}_{+}$, we have
\begin{equation*}
B_{0,q}=1,\quad\sum_{k=1}^{n}\binom{n}{k}_{q}\frac{1}{\lbrack k+1\rbrack_{q}}B_{n-k,q}=-B_{n,q},\,\,(n>0).
\end{equation*}
\end{prop}
By (\ref{eq:17}) and (\ref{eq:19}), we get
\begin{align}\label{eq:26}
p(x)&=\sum_{k=0}^{\infty}\frac{1}{\lbrack k\rbrack_{q}!}\left\langle\frac{e_{q}(t)-1}{t}t^{k}\Bigg\vert p(x)\right\rangle B_{k,q}(x)\\
&=\sum_{k=0}^{\infty}\frac{1}{\lbrack k\rbrack_{q}!}\left\langle\frac{e_{q}(t)-1}{t}\Bigg\vert t^{k}p(x)\right\rangle B_{k,q}(x)\nonumber\\
&=\sum_{k=0}^{\infty}\frac{1}{\lbrack k\rbrack_{q}!}B_{k,q}(x)\int_{0}^{1}t^{k}p(x)d_{q}x.\nonumber
\end{align}
It is known that
\begin{align}\label{eq:27}
\left(x-1\right)_{q}^{n}=\left(x-1\right)\left(x-q\right)\cdots (x-q^{n-1})\sim\left(e_{q}(t), t\right)
\end{align}
From (\ref{eq:17}) and (\ref{eq:27}), we have
\begin{align}\label{eq:28}
B_{n,q}(x)&=\sum_{k=0}^{n}\frac{1}{\lbrack k\rbrack_{q}!}\left\langle e_{q}(t)t^{k}\big\vert B_{n,q}(x)\right\rangle(x-1)_{q}^{k}\\
&=\sum_{k=0}^{n}\frac{1}{\lbrack k\rbrack_{q}!}\left\langle e_{q}(t)\big\vert t^{k}B_{n,q}(x)\right\rangle(x-1)_{q}^{k}\nonumber\\
&=\sum_{k=0}^{n}\binom{n}{k}_{q}B_{n-k,q}(1)(x-1)_{q}^{k}.\nonumber
\end{align}
From (\ref{eq:3}), we can derive
\begin{equation}\label{eq:29}
(x-1)_{q}^{n}=\sum_{m=0}^{n}\binom{n}{m}_{q}(-1)^{n-m}q^{\binom{n-m}{2}}x^{m}.
\end{equation}
Thus, by (\ref{eq:29}), we get
\begin{align}\label{eq:30}
t^{k}(x-1)_{q}^{n}&=\sum_{m=k}^{n}\binom{n}{m}_{q}(-1)^{n-m}q^{\binom{n-m}{2}}\frac{\lbrack m\rbrack_{q}!}{\lbrack m-k\rbrack_{q}!}x^{m-k}\\
&=\frac{\lbrack n\rbrack_{q}!}{\lbrack n-k\rbrack_{q}!}\sum_{m=0}^{n-k}\binom{n-k}{m}_{q}(-1)^{n-k-m}q^{\binom{n-k-m}{2}}x^{m}\nonumber\\
&=\frac{\lbrack n\rbrack_{q}!}{\lbrack n-k\rbrack_{q}!}(x-1)_{q}^{n-k}.\nonumber
\end{align}
By (\ref{eq:17}) and (\ref{eq:30}), we get
\begin{align}\label{eq:31}
(x-1)_{q}^{n}&=\sum_{k=0}^{n}\frac{1}{\lbrack k\rbrack_{q}!}\left\langle\frac{e_{q}(t)-1}{t}t^{k}\Bigg\vert(x-1)_{q}^{n}\right\rangle B_{k,q}(x)\\
&=\sum_{k=0}^{n}\binom{n}{k}_{q}B_{k,q}(x)\left\langle\frac{e_{q}(t)-1}{t}\Bigg\vert (x-1)_{q}^{n-k}\right\rangle\nonumber\\
&=\sum_{k=0}^{n}\binom{n}{k}_{q}B_{k,q}(x)\int_{0}^{1}(x-1)_{q}^{n-k}d_{q}x\nonumber\\
&=\sum_{k=0}^{n}\sum_{m=0}^{n-k}\binom{n}{k}_{q}\binom{n-k}{m}_{q}B_{k,q}(x)(-1)^{n-k-m}q^{\binom{n-k-m}{2}}\frac{1}{\lbrack m+1\rbrack_{q}}.\nonumber
\end{align}
From (\ref{eq:6}) and (\ref{eq:10}), we note that
\begin{equation}\label{eq:32}
\left\langle\frac{t}{e_{q}(t)-1}\Bigg\vert x^{n}\right\rangle=\sum_{k=0}^{\infty}\frac{B_{k,q}}{\lbrack k\rbrack_{q}!}\left\langle t^{k}\big\vert x^{n}\right\rangle=B_{n,q}.
\end{equation}
Let $\mathbb{P}_{n}=\left\{p(x)\in\mathbf{C}\lbrack x\rbrack\vert\deg{p(x)}\leq n\right\}$.\\
For $p(x)\in\mathbb{P}_{n}$, let us assume that
\begin{equation}\label{eq:33}
p(x)=\sum_{k=0}^{n}b_{k,q}B_{k,q}(x).
\end{equation}
By (\ref{eq:19}), we see that
\begin{equation}\label{eq:34}
\left\langle\left(\frac{e_{q}(t)-1}{t}\right)t^{k}\Bigg\vert B_{n,q}(x)\right\rangle=\lbrack n\rbrack_{q}!\delta_{n,k},\,\,(n,k\geq 0).
\end{equation}
Thus, from (\ref{eq:33}) and (\ref{eq:34}), we have
\begin{align}\label{eq:35}
\left\langle\left(\frac{e_{q}(t)-1}{t}\right)t^{k}\Bigg\vert p(x)\right\rangle&=\sum_{l=0}^{n}b_{l,q}\left\langle\left(\frac{e_{q}(t)-1}{t}\right)t^{k}\Bigg\vert B_{l,q}(x)\right\rangle\\
&=\sum_{l=0}^{n}b_{l,q}\lbrack l\rbrack_{q}!\delta_{l,k}=\lbrack k\rbrack_{q}!b_{k,q}.\nonumber
\end{align}
From (\ref{eq:16}), (\ref{eq:23}) and (\ref{eq:35}), we have
\begin{align}\label{eq:36}
b_{k,q}&=\frac{1}{\lbrack k\rbrack_{q}!}\left\langle\left(\frac{e_{q}(t)-1}{t}\right)t^{k}\Bigg\vert p(x)\right\rangle=\frac{1}{\lbrack k\rbrack_{q}!}\left\langle\frac{e_{q}(t)-1}{t}\Bigg\vert D_{q}^{k}p(x)\right\rangle\\
&=\frac{1}{\lbrack k\rbrack_{q}!}\int_{0}^{1}p^{(k)}(x)dx,\,\,\text{where}\,\,p^{(k)}(x)=D_{q}^{k}p(x).\nonumber
\end{align}
Therefore, by (\ref{eq:33}) and (\ref{eq:36}), we obtain the following theorem.

\begin{thm}\label{eq:thm2}
For $p(x)\in\mathbb{P}_{n}$, let $p(x)=\sum_{k=0}^{n}b_{k,q}B_{k,q}(x)$. Then we have
\begin{equation*}
b_{k,q}=\frac{1}{\lbrack k\rbrack_{q}!}\left\langle\frac{e_{q}(t)-1}{t}\Bigg\vert p^{(k)}(x)\right\rangle=\frac{1}{\lbrack k\rbrack_{q}!}\int_{0}^{1}p^{(k)}(x)d_{q}x,
\end{equation*}
where $p^{(k)}(x)=D_{q}^{k}p(x)$.
\end{thm}

\noindent Let us consider the $q$-Bernoulli polynomials of order $r$ as follows:
\begin{align}\label{eq:37}
\left(\frac{t}{e_{q}(t)-1}\right)^{r}e_{q}(xt)&=\underbrace{\left(\frac{t}{e_{q}(t)-1}\right)\times\cdots\times\left(\frac{t}{e_{q}(t)-1}\right)}_{r-\text{times}}e_{q}(xt)\\
&=\sum_{n=0}^{\infty}B_{n,q}^{(r)}(x)\frac{t^{n}}{\lbrack n\rbrack_{q}!}.\nonumber
\end{align}
In the special case, $x=0$, $B_{n,q}^{(r)}(0)=B_{n,q}^{(r)}$ is called the $n$-th $q$-Bernoulli number of order $r$. It is easy to show that
\begin{equation}\label{eq:38}
\left\langle\left(\frac{t}{e_{q}(t)-1}\right)^{r}\Bigg\vert x^{n}\right\rangle=\sum_{k=0}^{\infty}\frac{B_{k,q}^{(r)}}{\lbrack k\rbrack_{q}!}\left\langle t^{k}\big\vert x^{n}\right\rangle=B_{n,q}^{(r)}.
\end{equation}
From (\ref{eq:13}), (\ref{eq:32}) and (\ref{eq:38}), we note that
\begin{align}\label{eq:39}
B_{n,q}^{(r)}&=\left\langle\left(\frac{t}{e_{q}(t)-1}\right)^{r}\Bigg\vert x^{n}\right\rangle\\
&=\sum_{i_{1}+\cdots+i_{r}=n}^{}\binom{n}{i_{1},\cdots,i_{r}}_{q}\left\langle\frac{t}{e_{q}(t)-1}\Bigg\vert x^{i_{1}}\right\rangle\cdots\left\langle\frac{t}{e_{q}(t)-1}\Bigg\vert x^{i_{r}}\right\rangle\nonumber\\
&=\sum_{i_{1}+\cdots+i_{r}=n}^{}\binom{n}{i_{1},\cdots,i_{r}}_{q}B_{i_{1},q}\cdots B_{i_{r},q}.\nonumber
\end{align}
Therefore, by (\ref{eq:39}), we have the following lemma.
\begin{lem}\label{eq:lem3}
For $n\geq 0$, we have
\begin{equation*}
B_{n,q}^{(r)}=\sum_{i_{1}+\cdots+i_{r}=n}^{}\binom{n}{i_{1},\cdots,i_{r}}_{q}B_{i_{1},q}\cdots B_{i_{r},q}.
\end{equation*}
\end{lem}
By (\ref{eq:37}), we easily get
\begin{equation}\label{eq:40}
B_{n,q}^{(r)}(x)\sim\left(\left(\frac{t}{e_{q}(t)-1}\right)^{r}, t\right)
\end{equation}
and
\begin{equation}\label{eq:41}
B_{n,q}^{(r)}(x)=\left(\frac{t}{e_{q}(t)-1}\right)^{r}x^{n},\,\,\,\text{where}\,\, n,r\in\mathbf{Z}_{+}.
\end{equation}
Let us take $p(x)=B_{n,q}^{(r)}(x)=\sum_{k=0}^{n}\binom{n}{k}_{q}B_{n-k,q}^{(r)}x^{k}\in\mathbb{P}_{n}$. Then we may write
\begin{equation}\label{eq:42}
p(x)=B_{n,q}^{(r)}(x)=\sum_{k=0}^{n}b_{k,q}B_{k,q}(x).
\end{equation}
From (\ref{eq:42}), we have
\begin{align}\label{eq:43}
p^{(k)}(x)=D_{q}^{k}B_{n,q}^{(r)}(x)&=\lbrack n\rbrack_{q}\lbrack n-1\rbrack_{q}\cdots\lbrack n-k+1\rbrack_{q}B_{n-k,q}^{(r)}(x)\\
&=\lbrack k\rbrack_{q}!\binom{n}{k}_{q}B_{n-k,q}^{(r)}(x).\nonumber
\end{align}
By (\ref{eq:36}) and (\ref{eq:43}), we get
\begin{align}\label{eq:44}
b_{k,q}&=\frac{1}{\lbrack k\rbrack_{q}!}\left\langle\left(\frac{e_{q}(t)-1}{t}\right)t^{k}\Bigg\vert p(x)\right\rangle=\frac{1}{\lbrack k\rbrack_{q}!}\left\langle\frac{e_{q}(t)-1}{t}\Bigg\vert D_{q}^{k}p(x)\right\rangle\\
&=\binom{n}{k}_{q}\left\langle\frac{e_{q}(t)-1}{t}\Bigg\vert B_{n-k,q}^{(r)}(x)\right\rangle=\binom{n}{k}_{q}\left\langle t^{0}\Bigg\vert\left(\frac{t}{e_{q}(t)-1}\right)^{r-1}x^{n-k}\right\rangle\nonumber\\
&=\binom{n}{k}_{q}B_{n-k,q}^{(r-1)}\nonumber.
\end{align}
Therefore, by Theorem \ref{eq:thm2} and (\ref{eq:42}), we obtain the following theorem.

\begin{thm}\label{eq:thm4}
For $n\geq 0$, we have
\begin{align*}
B_{n,q}^{(r)}(x)&=\sum_{k=0}^{n}\binom{n}{k}_{q}\left\langle\frac{e_{q}(t)-1}{t}\Bigg\vert B_{n-k,q}^{(r)}(x)\right\rangle B_{k,q}(x)\\
&=\sum_{k=0}^{n}\binom{n}{k}_{q}B_{n-k,q}^{(r-1)}B_{k,q}(x).
\end{align*}
\end{thm}

\noindent For $p(x)\in\mathbb{P}_{n}$, let us assume that
\begin{equation}\label{eq:45}
p(x)=\sum_{k=0}^{n}b_{k,q}^{(r)}B_{k,q}^{(r)}(x).
\end{equation}
By (\ref{eq:40}), we easily get
\begin{equation}\label{eq:46}
\left\langle\left(\frac{e_{q}(t)-1}{t}\right)^{r}t^{k}\Bigg\vert B_{n,q}^{(r)}(x)\right\rangle=\lbrack n\rbrack_{q}!\delta_{n,k},\,\,\,(n,k\geq 0).
\end{equation}
From (\ref{eq:45}) and (\ref{eq:46}), we have
\begin{align}\label{eq:47}
\left\langle\left(\frac{e_{q}(t)-1}{t}\right)^{r}t^{k}\Bigg\vert p(x)\right\rangle&=\sum_{l=0}
^{n}b_{l,q}^{(r)}\left\langle\left(\frac{e_{q}(t)-1}{t}\right)^{r}t^{k}\Bigg\vert B_{l,q}^{(r)}(x)\right\rangle\\
&=\sum_{l=0}^{n}b_{l,q}^{(r)}\lbrack l\rbrack_{q}!\delta_{l,k}=\lbrack k\rbrack_{q}!b_{k,q}^{(r)}.\nonumber
\end{align}
By (\ref{eq:47}), we get
\begin{equation}\label{eq:48}
b_{k,q}^{(r)}=\frac{1}{\lbrack k\rbrack_{q}!}\left\langle\left(\frac{e_{q}(t)-1}{t}\right)^{r}t^{k}\Bigg\vert p(x)\right\rangle.
\end{equation}
Therefore, by (\ref{eq:45}) and (\ref{eq:48}), we obtain the following theorem.

\begin{thm}\label{eq:thm5}
For $p(x)\in\mathbb{P}_{n}$, let $p(x)=\sum_{k=0}^{n}b_{k,q}^{(r)}B_{k,q}^{(r)}(x)$. Then we have
\begin{align*}
b_{k,q}^{(r)}=\frac{1}{\lbrack k\rbrack_{q}!}\left\langle\left(\frac{e_{q}(t)-1}{t}\right)^{r}t^{k}\Bigg\vert p(x)\right\rangle.
\end{align*}
\end{thm}

\noindent Let us take $p(x)=B_{n,q}(x)$. Then, by Theorem \ref{eq:5}, we get
\begin{equation}\label{eq:49}
B_{n,q}(x)=p(x)=\sum_{k=0}^{n}b_{k,q}^{(r)}B_{k,q}^{(r)}(x),
\end{equation}
where
\begin{equation}\label{eq:50}
b_{k,q}^{(r)}=\frac{1}{\lbrack k\rbrack_{q}!}\left\langle\left(\frac{e_{q}(t)-1}{t}\right)^{r}t^{k}\Bigg\vert p(x)\right\rangle=\frac{1}{\lbrack k\rbrack_{q}!}\left\langle\left(\frac{e_{q}(t)-1}{t}\right)^{r}t^{k}\Bigg\vert B_{n,q}(x)\right\rangle.
\end{equation}
For $k<r$, by (\ref{eq:50}), we have
\begin{align}\label{eq:51}
b_{k,q}^{(r)}&=\frac{1}{\lbrack k\rbrack_{q}!}\left\langle\left(e_{q}(t)-1\right)^{r}\frac{1}{t^{r-k}}\Bigg\vert B_{n,q}(x)\right\rangle\\
&=\frac{1}{\lbrack k\rbrack_{q}!}\left(\frac{1}{\lbrack n+r-k\rbrack_{q}\cdots\lbrack n+1\rbrack_{q}}\right)\left\langle\left(e_{q}(t)-1\right)^{r}\left(\frac{1}{t}\right)^{r-k}\Bigg\vert t^{r-k}B_{n+r-k,q}(x)\right\rangle\nonumber\\
&=\left(\frac{1}{\lbrack k\rbrack_{q}!\lbrack r-k\rbrack_{q}!}\right)\left(\frac{\lbrack r-k\rbrack_{q}!}{\lbrack n+r-k\rbrack_{q}\cdots\lbrack n+1\rbrack_{q}}\right)\left\langle\left(e_{q}(t)-1\right)^{r}\big\vert B_{n+r-k,q}(x)\right\rangle\nonumber\\
&=\frac{1}{\lbrack r\rbrack_{q}!}\frac{\binom{r}{k}_{q}}{\binom{n+r-k}{r-k}_q}\sum_{j=0}^{r}\binom{r}{j}(-1)^{r-j}\left\langle\left(e_{q}(t)\right)^{j}\big\vert B_{n+r-k,q}(x)\right\rangle\nonumber
\end{align}
\begin{align*}
&=\frac{1}{\lbrack r\rbrack_{q}!}\frac{\binom{r}{k}_{q}}{\binom{n+r-k}{r-k}_q}\sum_{j=0}^{r}\binom{r}{j}(-1)^{r-j}\sum_{m=0}^{n+r-k}\sum_{m_{1}+\cdots+m_{j}=m}^{}\binom{m}{m_{1},\cdots ,m_{j}}_{q}\\
&\quad\times\binom{n+r-k}{m}_{q}B_{n+r-k-m,q}.
\end{align*}
Let us assume that $k\geq r$. Then, by (\ref{eq:50}), we get
\begin{align}\label{eq:52}
b_{k,q}^{(r)}&=\frac{1}{\lbrack k\rbrack_{q}!}\left\langle\left(e_{q}(t)-1\right)^{r}\big\vert t^{k-r}B_{n,q}(x)\right\rangle\\
&=\frac{1}{\lbrack k\rbrack_{q}!}\lbrack n\rbrack_{q}\lbrack n-1\rbrack_{q}\cdots \lbrack n-k+r+1\rbrack_{q}\left\langle\left(e_{q}(t)-1\right)^{r}\big\vert B_{n-k+r,q}(x)\right\rangle\nonumber\\
&=\frac{\lbrack k-r\rbrack_{q}!}{\lbrack k\rbrack_{q}!}\binom{n}{k-r}_{q}\sum_{j=0}^{r}\binom{r}{j}(-1)^{r-j}\left\langle\left(e_{q}(t)\right)^{j}\big\vert B_{n-k+r,q}(x)\right\rangle\nonumber\\
&=\frac{1}{\lbrack r\rbrack_{q}!}\frac{\binom{n}{k-r}_{q}}{\binom{k}{r}_{q}}\sum_{j=0}^{r}\binom{r}{j}(-1)^{r-j}\sum_{m=0}^{n-k+r}\sum_{m_{1}+\cdots+m_{j}=m}^{}\binom{m}{m_{1},\cdots ,m_{j}}_{q}\nonumber\\
&\quad\times\frac{\left\langle t^{m}\big\vert B_{n-k+r,q}(x)\right\rangle}{\lbrack m\rbrack_{q}!}\nonumber\\
&=\frac{1}{\lbrack r\rbrack_{q}!}\frac{\binom{n}{k-r}_{q}}{\binom{k}{r}_{q}}\sum_{j=0}^{r}\binom{r}{j}(-1)^{r-j}\sum_{m=0}^{n-k+r}\sum_{m_{1}+\cdots+m_{j}=m}^{}\binom{m}{m_{1},\cdots ,m_{j}}_{q}\nonumber\\
&\quad\times\binom{n-k+r}{m}_{q}B_{n-k+r-m,q}.\nonumber
\end{align}
Therefore, by (\ref{eq:49}), (\ref{eq:51}) and (\ref{eq:52}), we obtain the following theorem.

\begin{thm}\label{eq:thm6}
For $n\in\mathbf{Z}_{+}$ and $r\in\mathbf{N}$, we have
\begin{align*}
B_{n,q}(x)&=\sum_{k=0}^{r-1}\frac{1}{\lbrack r\rbrack_{q}!}\frac{\binom{r}{k}_{q}}{\binom{n+r-k}{r-k}_{q}}\Bigg\{\sum_{j=0}^{r}\binom{r}{j}(-1)^{r-j}\sum_{m=0}^{n-k+r}\sum_{m_{1}+\cdots+m_{j}=m}^{}\binom{m}{m_{1},\cdots ,m_{j}}_{q}\\
&\quad\times\binom{n-k+r}{m}_{q}B_{n+r-k-m,q}\Bigg\}B_{k,q}^{(r)}(x)+\sum_{k=r}^{n}\frac{\binom{n}{k-r}_{q}}{\lbrack r\rbrack_{q}!\binom{r}{k}_{q}}\\
&\quad\times\Bigg\{\sum_{j=0}^{r}\binom{r}{j}(-1)^{r-j}\sum_{m=0}^{n-k+r}\sum_{m_{1}+\cdots+m_{j}+m}^{}\binom{m}{m_{1},\cdots ,m_{j}}_{q}\binom{n-k+r}{m}_{q}\\
&\quad\times B_{n-k+r-m,q}\Bigg\}B_{k,q}^{(r)}(x).
\end{align*}
\end{thm}


\noindent
\author{Department of Mathematics, Sogang University, Seoul 121-742, Republic of Korea
\\e-mail: dskim@sogang.ac.kr}\\
\\
\author{Department of Mathematics, Kwangwoon University, Seoul 139-701, Republic of Korea
\\e-mail: tkkim@kw.ac.kr}
\end{document}